\documentclass[12pt]{article}
\usepackage{amssymb,amsmath,latexsym}

\DeclareMathOperator{\aut}{Aut}

\DeclareMathOperator{\Mat}{Mat} \DeclareMathOperator{\orb}{Orb}
 
 \DeclareMathOperator{\sym}{Sym}

\def\C{{\mathbb C}}

\def\N{{\mathbb N}}

\def\RR{{\mathbb R}}

\def\A{{\cal A}}

\def\D{{\cal D}}

\def\R{{\cal R}}

\def\proof{{\bf Proof}.\ }
\def\bull{\vrule height .9ex width .8ex depth -.1ex }

\makeatletter
\renewcommand{\subsection}{\@startsection{subsection}{2}{0mm}{-2mm}{-2mm}
{\bf\normalsize}}
\def\sbsn{\subsection{\hspace{-3mm}}}
\def\sbst#1{\subsection*{#1}}

\makeatother

\newtheorem{formula}{}[section]
\newtheorem{proposition}[formula]{Proposition}
\newtheorem{definition}[formula]{Definition}
\newtheorem{corollary}[formula]{Corollary}
\newtheorem{remark}[formula]{Remark}
\newtheorem{lemma}[formula]{Lemma}
\newtheorem{theorem}[formula]{Theorem}

\newtheorem{hypothesis}[formula]{Conjecture}
\newtheorem{example}[formula]{Example}

\def\thrm{\begin{theorem}}
\def\thrml#1{\begin{theorem}\label{#1}}
\def\ethrm{\end{theorem}}
\def\prpstn{\begin{proposition}}
\def\prpstnl#1{\begin{proposition}\label{#1}}
\def\eprpstn{\end{proposition}}
\def\rmrk{\begin{remark}}
\def\rmrkl#1{\begin{remark}\label{#1}}
\def\ermrk{\end{remark}}
\def\dfntn{\begin{definition}}
\def\dfntnl#1{\begin{definition}\label{#1}}
\def\edfntn{\end{definition}}
\def\nmrt{\begin{enumerate}}
\def\enmrt{\end{enumerate}}
\def\tm#1{\item[{\rm (#1)}]}
\def\qtn{\begin{equation}}
\def\qtnl#1{\begin{equation}\label{#1}}
\def\eqtn{\end{equation}}
\def\lmm{\begin{lemma}}
\def\lmml#1{\begin{lemma}\label{#1}}
\def\elmm{\end{lemma}}
\def\crllr{\begin{corollary}}
\def\crllrl#1{\begin{corollary}\label{#1}}
\def\ecrllr{\end{corollary}}
\def\hpthss{\begin{hypothesis}}
\def\hpthssl#1{\begin{hypothesis}\label{#1}}
\def\ehpthss{\end{hypothesis}}
\def\xmpl{\begin{example}}
\def\xmpll#1{\begin{example}\label{#1}}
\def\exmpl{\end{example}}
\def\css{\begin{cases}}
\def\ecss{\end{cases}}

\begin{document}

\begin{center}
{\large \bf {On amorphic $C$-algebras }}

\smallskip

\medskip\medskip
\end{center}

\textbf{I.N. Ponomarenko \footnote{Partially supported by RFFI,
grants 03-01-000349, 05-01-00-899,  NSH-2251.2003.1. The author
thanks the Institute for Advanced Studies in Basic Sciences
(IASBS) for its hospitality during the time that this paper was
written}}

Petersburg Department of V.A. Steklov

Institute of Mathematics

Fontanka 27, St. Petersburg 191023, Russia

inp@pdmi.ras.ru

http://www.pdmi.ras.ru/\~{}inp

\smallskip\medskip

\textbf{A. Rahnamai Barghi}

Institute for Advanced Studies in Basic Sciences

Gava Zang, P.O.Box: 45195-1159 Zanjan, Iran

\qquad\qquad\qquad

email: rahnama@iasbs.ac.ir

\smallskip\newpage

running title: On amorphic $C$-algebras

\smallskip\medskip

SEND PROOFS TO:\medskip

\textbf{A. Rahnamai Barghi}

Institute for Advanced Studies in Basic Sciences

Gava Zang, P.O.Box: 45195-1159 Zanjan, Iran

email: rahnama@iasbs.ac.ir

\newpage

{\large ABSTRACT}

\medskip

An amorphic association scheme has the property that any of its
fusion is also an association scheme. In this paper we generalize
the property to be amorphic to an arbitrary $C$-algebra and prove
that any amorphic $C$-algebra is determined up to isomorphism by
the multiset of its degrees and an additional integer equal $\pm
1$. Moreover, we show that any amorphic $C$-algebra with rational
structure constants is the fusion of an amorphic homogeneous
$C$-algebra. Finally, as a special case of our results we obtain
the well-known Ivanov's characterization of intersection numbers
of amorphic association schemes.\\[2mm]

\medskip

\noindent{\it AMS Subject Classification 2000: 20C05 05E30}

\section{Introduction}
One of the important problems in association schemes theory
consists of the constructive enumeration of all the subschemes of
a given scheme. Such a problem becomes intractable when the number
of subschemes is large. In the limit case this number coincides
with the number of all partitions of the set of classes of the
scheme in question. In order to describe this situation, the
notion of amorphic cellular ring (the adjacency ring of an
association scheme) was introduced in~\cite{GK85}. In paper
\cite{I85} that is a continuation of~\cite{GK85}, a complete
characterization of the intersection numbers array of an arbitrary
amorphic cellular ring is given. Since then a number of examples
of amorphic association schemes were found (see~\cite{DX,IMY91}).
However, the problem of characterization of amorphic association
schemes is far from being solved.

The notion of $C$-algebra goes back to the old paper of Kawada
(1942) in which a finite dimensional analog of the Tannaka duality
was introduced (see Section \ref{l2} in the below for a background
of $C$-algebras). The two basic examples of $C$-algebras (which
are also table algebras in sense of \cite{AB}) are the algebra of
conjugacy classes of a finite group and the character algebra of a
finite group. Besides, in \cite{BI} it was showed that the
adjacency algebra of a commutative association scheme is a
$C$-algebra. A special case of non-commutative $C$-algebras
(including all the above examples and also the adjacency algebra
of an arbitrary association scheme) was studied in~\cite{AFM}.\\

The main purpose of the paper is to generalize the concept of
amorphic association scheme for $C$-algebras and study amorphic
$C$-algebras in details. According to our definition (see
Definition~\ref{l1}), the adjacency algebra of an amorphic
association scheme is an amorphic $C$-algebra, but not all of the
latter arise in this way. For instance, affine $C$-algebras
corresponding to affine planes of order $q$ (see Section~\ref{l2})
come from affine association schemes and so are amorphic due to
\cite{GK85}. On the other hand, there are some generalized affine
$C$-algebras (see Subsection \ref{l3}) which do not come from
association schemes.

After studying some elementary properties of amorphic
$C$-algebras, we prove that each such algebra is determined up to
isomorphism by the multiset of its degrees and an additional
integer equal $\pm 1$ (Theorem~\ref{f060205b}). As a special case
of this result we obtain the well-known Ivanov's characterization
of intersection numbers of amorphic association schemes
\cite{I85}. In Section \ref{l4}, the class of homogeneous amorphic
$C$-algebras is studied (this class contains generalized affine
$C$-algebras mentioned in the above). We show that any nontrivial
amorphic $C$-algebra with integer structure constants is a fusion
of a generalized affine algebra (Corollary \ref{l5}). Furthermore,
we give a necessary and sufficient condition for a $C$-algebra to be
homogeneous and amorphic (Theorem \ref{f150502b}). Finally it is
proved that the dual of a homogeneous amorphic table algebras is a
homogeneous amorphic table algebra too (Theorem \ref{f160205a}).

Concerning associaton schemes and permutation groups we refer the
reader to \cite{BI} and \cite{DIX} respectively.

\section{$C$-algebras}\label{l2}
In the definition below we follow the paper~\cite{EPV} in which
arbitrary non-commutative $C$-algebras were studied from the
duality theory point of view.

Let $A$ be an associative algebra over $\C$ with the unit $e$,
semilinear involution anti-automorphism $*$ and a linear basis
$\R$. The structure constants of this algebra with respect to the
basis $\R$ are defined by \qtnl{f090205a} r\cdot s=\sum_{t\in\R
}c_{r,s}^t t,\quad r,s\in\R. \eqtn

\dfntn The pair $\A=(A,\R)$ is called a $C$-algebra, if the
following conditions are satisfied: \nmrt \tm{C1} $e\in\R$ and
$r^*\in\R$ for all $r\in\R$, \tm{C2} $c_{r,s}^t\in\RR$ for all
$r,s,t\in\R$, \tm{C3} for all $r,s\in\R$ the equality
$c_{r,s}^e=d_r\delta_{r^*,s}$ holds where $d_r>0$ and $\delta$ is
the Kronecker symbol, \tm{C4} a linear functional on $A$, defined
by the mapping $r\mapsto d_r$, $r\in\R$, is a one-dimensional
$*$-representation of the algebra~$A$. \enmrt \edfntn The
$C$-algebra $\A$ is called {\it symmetric} (resp. {\it
commutative}) if the involution $*$ is trivial (resp. the algebra
$\A$ is commutative). Obviously, any symmetric $C$-algebra is
commutative.

Two $C$-algebras $\A=(A,\R)$ and~$\A=(A',\R')$ are called {\it
isomorphic}, if there exists an $*$-algebra isomorphism
$\varphi:A\to A'$ such that $\varphi(\R)=\R'$; such an isomorphism
is called the {\it isomorphism} from~$\A$ to~$\A'$. For $\A=\A'$
the set of all of them form a group $\aut(\A)$ called the {\it
automorphism group} of the $C$-algebra~$\A$. From the definition
it follows that this group acts faithfully on the set $\R$ and
hence it can be naturally identified with a subgroup of the group
$\sym(\R^\#)$ where $\R^\#=\R\setminus\{e\}$.

\sbst{Example: affine $C$-algebras (cfg. \cite{GK85}).} {\it An
affine plane} is a triple consisting of a set $P$ of points, a set
$L$ of lines and an incidence relation between the points and the
lines such that \nmrt \tm{A1} any two different points are
incident to a unique line, \tm{A2} for any point $p$ and for any
line $l$ not incident to $p$ there exists a unique line incident
to $p$ and parallel to $l$, \tm{A3} there exist three different
points not incident to the same line. \enmrt From now on we assume
that the plane is finite, i.e. the set $P\cup L$ is finite. Then
there exists an integer $q\ge 2$, called the {\it order} of the
plane, such that $|P|=q^2$ and $|L|=q^2+q$, and each line incident
to exactly $q+1$ points. Besides, there are exactly $q+1$ classes
$B_1,\ldots,B_{q+1}$ of pairwise parallel (non-intersecting)
lines, each of which is of cardinality~$q$.

Let $\Mat_P$ be the full matrix algebra of degree $|P|$ over $\C$
whose rows and columns are indexed by the
elements of $P$. Given $i=1,\ldots,q+1$ we define a \{0,1\}-matrix
$r_i$ belongs to $\Mat_P$ by
$$
(r_i)_{u,v}= \css
1, &\text{if $u\ne v$ and $l(u,v)\in B_i$},\\
0, &\text{otherwise}, \ecss
$$
where $l(u,v)$ is the line incident to both $u$ and $v$. One can see
that
\qtnl{f050205d}
r_ir_j= \css
(q-1)r_0+(q-2)r_i, &\text{if $i=j$},\\
\displaystyle\sum_{k\ne 0,i,j}r_k, &\text{if $i\ne j$}
\ecss
\eqtn
where $r_0$ is the identity matrix of $\Mat_P$. So the set
$\R=\{r_i:\ i=0,\ldots,q+1\}$ is a linear base of the subalgebra
$A$ of the algebra $\Mat_P$ generated by $\R$. Since conditions
(C1)--(C4) are trivially satisfied (with $*$ being the Hermitian
conjugation in $\Mat_P$), the pair $\A=(A,\R)$ is a $C$-algebra.
We call it the {\it affine} $C$-algebra. From (\ref{f050205d}),
it follows that any affine $C$-algebra is determined up to
isomorphism by the
order of the corresponding affine plane (and does not depend on
the choice of the plane). Moreover, it is easy to see that $\A$ is
symmetric and  $\aut(\A) \simeq \sym(q+1)$.\bull \vspace{2mm}

Let $\A=(A,\R)$ be a $C$-algebra. Then the structure constants
$c_{r,s}^t$ of the algebra $A$ with respect to the base $\R$
satisfy some relations. For instance, from (C4) it follows that
$d_e=1$ and $d_r=d_{r^*}$ for all $r\in\R$. Besides, the
associativity of $A$ exactly means that given $r,s,t\in\R$ the
following equalities hold: \qtnl{f050205a}
\sum_{v\in\R}c_{r,s}^vc_{v,t}^u=
\sum_{v\in\R}c_{r,v}^uc_{s,t}^v,\quad u\in\R. \eqtn In fact, the
left-hand side and the right-hand side of the equalities
(\ref{f050205a}) are the coefficients of $u$ in the products
$(rs)t$ and $r(st)$ respectively. Take $u=e$. Then from condition
(C3) it follows that \qtnl{f050205b}
d_tc_{r,s}^{t^*}=d_rc_{s,t}^{r^*}=d_sc_{t,r}^{s^*},\quad
r,s,t\in\R. \eqtn Finally, applying the linear functional from
condition (C4) to all the products of two elements of $\R$, we get
the equalities \qtnl{f050205c}
d_rd_s=\sum_{t\in\R}d_tc_{r,s}^{t},\quad r,s\in\R. \eqtn Below we
denote by $D(\A)$ the multiset $\{d_r:\ r\in\R\}$ and call the
number $n(\A)=\sum_{r\in\R}d_r$ the {\it degree} of $\A$. In
particular, if $\A$ is the affine $C$-algebra corresponding to an
affine plane of order $q$, then from (\ref{f050205d}) it follows
that $D(\A)=\{1,q-1,\ldots,q-1\}$ ($q+1$ times) and $n(\A)=q^2$.
This gives an example of a {\it homogeneous} $C$-algebra, i.e. one
for which the number $d_r$ does not depend on $r\in\R^\#$.

Let $\pi$ be a partition of $\R$ such that $\{e\}\in\pi$.
Then the pair $\A_\pi=(A_\pi,\R_\pi)$ where $A_\pi$ is the linear
subspace of $A$ spanned by the set
$$
\R_\pi=\{\sum_{r\in R}r:\ R\in\pi\},
$$
is called the {\it $\pi$-fusion} of $\A$. It is easy to see that
$\A_\pi$ is a $C$-algebra iff $A_\pi$ is an $*$-subalgebra of $A$.
For example, the latter is always true whenever
$|\pi|\in\{2,|\R|\}$. In this paper we are interesting in the case
when each fusion of $\A$ is a $C$-algebra.

\section{A characterization of amorphic $C$-algebras}

\sbsn The following definition extends the notion of amorphic
association scheme to arbitrary $C$-algebras.

\dfntnl{l1}
A $C$-algebra $\A=(A,\R)$ is called amorphic if for
each partition $\pi$ of $\R$ the $\pi$-fusion
$\A_\pi=(A_\pi,\R_\pi)$ of $\A$ is a $C$-algebra.
\edfntn

From the definition it follows that the adjacency algebra of an
amorphic association scheme is an amorphic $C$-algebra. Next,
obviously, any $C$-algebra of dimension at most $3$ is amorphic;
we call them {\it trivial} amorphic $C$-algebras. Nontrivial
examples of amorphic $C$-algebras can be obtained from
\cite[Theorem~3.3]{GK85} which implies that any affine $C$-algebra
is amorphic (see also paper \cite{DX} and the references in it). A
full characterization of amorphic $C$-algebras isomorphic to the
adjacency algebras of association schemes follows from~\cite{I85}
(a full characterization of amorphic cyclotomic association
schemes over finite fields up to association scheme isomorphisms
follows from~\cite{BMW}). However, there is a number of $C$-algebras
which are not isomorphic to the adjacency algebra of an
association scheme (e.g. the structure constants of the latter
must be nonnegative integers) and our results show that such
examples can be found in the class of amorphic $C$-algebras.

Let $\A=(A,\R)$ be an amorphic $C$-algebra. Given $r,s\in\R^\#$
set $\pi=\pi_{r,s}$ to be the partition of $\R$ with classes
$\{e\}$, $\{r\}$, $\{s\}$ and $\R^\#\setminus\{e,r,s\}$ (we
don't assume that $r\ne s$). Since
$\A_\pi=(A_\pi,\R_\pi)$ is a $C$-algebra with $r,s\in\R_\pi$, we
see that \qtnl{f050205i} c_{r,s}^u=c_{r,s}^v,\quad
u,v\in\R^\#\setminus\{r,s\}. \eqtn Moreover, for $r=s$ we also
have $r^*\in\{r,t\}$ where $t$ is the sum of all
$u\in\R^\#\setminus\{r\}$. If $|\R|\ge 4$, this implies by
condition (C1) that $r=r^*$. Thus any nontrivial amorphic
$C$-algebra is symmetric and hence any amorphic $C$-algebra is
commutative (c.f. \cite[Theorem~2.3]{GK85}).

\thrml{f050205g}
Let $\A=(A,\R)$ be a symmetric $C$-algebra. Then
the following three statements are equivalent:
\nmrt \tm{1} $\A$ is amorphic,
\tm{2} given $r,s\in\R^\#$ the $\pi_{r,s}$-fusion of
$\A$ is a $C$-algebra,
\tm{3} given $r,s\in\R^\#$ the equalities (\ref{f050205i}) hold.
\enmrt
\ethrm
\proof  The implications
$(1)\Rightarrow(2)$ and $(2)\Rightarrow(3)$ are obvious. To prove
the implication $(3)\Rightarrow(1)$ let $\pi$ be a partition of
$\R$ with $\{e\}\in\pi$. Then the linear space $A_\pi$ is closed
with respect to $*$ because the $C$-algebra $\A$ is symmetric.
Given $R,S\in\pi$ set $T=\R\setminus (R\cup S\cup\{e\})$.
Then from (\ref{f050205i}) it follows that
$$
\sum_{r\in R}\sum_{s\in S}c_{r,s}^u= \sum_{r\in R}\sum_{s\in
S}c_{r,s}^v,\quad u,v\in T.
$$
This implies that $\sum_{r\in R}r\sum_{s\in S}s\in A_\pi$ and so
the linear space $A_\pi$ is an $*$-subalgebra of $A$. Thus the
fusion $\A_\pi$ is a $C$-algebra and we are done.\bull

\sbsn\label{f090205b} The following theorem gives explicit
formulas for the structure constants of a nontrivial amorphic
$C$-algebra $(A,\R)$. We observe that due to equalities
(\ref{f050205b}) it suffices to find the structure constants of
the form $c_{r,s}^t$ where $r,s,t\in\R^\#$ are such that $r\ne
s\ne t\ne r$ or $r=s\ne t$ or $r=s=t$.

\thrml{f060205b}
Let $\A=(A,\R)$ be a nontrivial amorphic $C$-algebra. Then the
structure constants of $\A$ are uniquely determined by the multiset
$D(\A)$ and by the number $\varepsilon\in\{-1,+1\}$. More exactly,
given $r,s,t\in\R^\#$ we have
\qtnl{f070205a}
c_{r,s}^t=
\css
\frac{d_rd_s}{(\sqrt{n}+\varepsilon)^2}, &\text{if $r\ne s\ne t\ne r$},\\[4pt]
\phantom{-}\frac{d_r}{2}\left(1+\frac{d_r-d'_r}{(\sqrt{n}+\varepsilon)^2}\right),
&\text{if $r=s\ne t$},\\[4pt]
-\frac{d'_r}{2}\left(1+\frac{d_r-d'_r}{(\sqrt{n}+\varepsilon)^2}\right)
+d_r-1, &\text{if $r=s=t$},\\
\ecss
\eqtn
where $n=n(\A)$ and $d'_r=n-1-d_r$.
\ethrm
\proof Given two distinct elements $r,s\in\R^\#$ set $\pi=\pi_{r,s}$ and
$t$ to be the sum of all $u\in\R^\#\setminus\{r,s\}$. Then
$\A_\pi=(A_\pi,\R_\pi)$ is an amorphic $C$-algebra of dimension
$4$ with $\R_\pi=\{e,r,s,t\}$. Since it is symmetric, formulas
(\ref{f050205c}) and (\ref{f050205i}) imply that
\qtnl{f070205b}
\left\{\begin{array}{rcl}
d_r^2-d_r & = & d_rc_{r,r}^r+d_sc_{r,r}^s+d_tc_{r,r}^t\\
d_rd_s    & = & d_rc_{r,s}^r+d_sc_{r,s}^s+d_tc_{r,s}^t\\
0         & = & d_tc_{r,r}^s-d_sc_{r,r}^t\\
\end{array}\right.
\eqtn
where $d_t=n-1-d_r-d_s$. By (\ref{f050205b}) the number
$d_wc_{u,v}^w$ does not depend on $u,v,w\in\R_\pi^\#$ such that
$\{u,v,w\}=\{r,s,t\}$. Denote it by $T$ and set
$X_{uv}=d_vc_{u,u}^v$ where $u,v\in\{r,s,t\}$. Then equations
(\ref{f070205b}) produce the following system of linear equations:
$$
\begin{pmatrix}
 1   & 0   & 0   & 1   &  1   & 0   &  0   & 0   &  0\\
 0   & 1   & 0   & 0   &  0   & 1   &  1   & 0   &  0\\
 0   & 0   & 1   & 0   &  0   & 0   &  0   & 1   &  1\\
 0   & 0   & 0   & 1   &  0   & 0   &  1   & 0   &  0\\
 0   & 0   & 0   & 0   &  1   & 0   &  0   & 1   &  0\\
 0   & 0   & 0   & 0   &  0   & 1   &  0   & 0   &  1\\
 0   & 0   & 0   & d_t & -d_s & 0   &  0   & 0   &  0\\
 0   & 0   & 0   & 0   &  0   & d_r & -d_t & 0   &  0\\
 0   & 0   & 0   & 0   &  0   & 0   &  0   & d_s & -d_r
\end{pmatrix}
\begin{pmatrix}
X_{rr}\\
X_{ss}\\
X_{tt}\\
X_{rs}\\
X_{rt}\\
X_{st}\\
X_{sr}\\
X_{tr}\\
X_{ts}
\end{pmatrix}
=
\begin{pmatrix}
d_r^2-d_r\\
d_s^2-d_s\\
d_t^2-d_t\\
d_rd_s-T\\
d_td_r-T\\
d_sd_t-T\\
0\\
0\\
0
\end{pmatrix}
$$
This system has the unique solution from which we obtain the
following expressions for the structure constants:
\qtnl{f070205d}
c_{u,u}^v=c_{u,u}^t=T\frac{d_u-d'_u}{2d_td_v}+\frac{d_u}{2},\ \ \
c_{u,u}^u=T\frac{d'_u(d'_u-d_u)}{2m}+\frac{2d_u-2-d'_u}{2},
\eqtn
where $\{u,v\}=\{r,s\}$ and $m=d_rd_sd_t$. By substituting
these solutions to equality (\ref{f050205a}) with $u=r$ we obtain
the following quadratic equation with respect to $T$:
$$
(n-1)^2T^2-2m(n+1)T+m^2=0.
$$
Thus $T=m/(\sqrt{n}+\varepsilon_{r,s})^2$ where
$\varepsilon_{r,s}=\pm 1$. (We note that due to condition (C3) we
have $\sqrt{n}+\varepsilon_{r,s}>0$.) Taking $u=r,s$ in the second
of equalities (\ref{f070205d}) we see that in fact the number
$\epsilon_{r,s}$ does not depend on the choice of~$r,s$. Denote it
by $\varepsilon$. Then by the definition of $T$ this proves the
first equality in (\ref{f070205a}); the other ones follow
from~(\ref{f070205d}) because due to Theorem~\ref{f050205g} we
have $c_{r,s}^w=c_{r,s}^t$ for all $w\in\R^\#\setminus\{r,s\}$.\bull \\

Let us rewrite formulas (\ref{f070205a}) in slightly different
form. For $r\in\R$ set $g_r=m+\varepsilon$ where $m=\sqrt{n}$.
Then a straightforward computation shows that \qtnl{f200205a}
d_r=g_r(m+\varepsilon),\quad
c_{r,r}^r=(g_r+\varepsilon)(g_r+2\varepsilon)-\varepsilon
m-2,\quad c_{r,r}^s=g_r(g_r+\varepsilon) \eqtn for all
$s\in\R\setminus\{r\}$. Suppose from now that the $C$-algebra $\A$
is the adjacency algebra of an association scheme. Then it's
structure constants are nonnegative integers. This implies that
the number $n$ is an exact square and relations (\ref{f200205a})
show that each basis graph of the association scheme is a strongly
regular graph of Latin square type (for $\varepsilon=-1$) or of
negative Latin square type (for $\varepsilon=1$). Thus as a
special case of Theorem~\ref{f060205b} we obtain the result of
A.~V.~Ivanov  who proved in~\cite{I85} that the number of points
of a nontrivial amorphic association scheme is an exact square and
all basis graphs of the scheme are simultaneously strongly regular
graphs of Latin square type or of negative Latin square type.

Theorem~\ref{f060205b} enables us to compute the automorphism
group of an arbitrary amorphic $C$-algebra $\A=(A,\R)$. To do this
denote by $\aut(D(\A))$ the group of all permutations
$g\in\sym(\R)$ such that $e^g=e$ and $d_r=d_{r^g}$ for all
$r\in\R^\#$. Clearly,
$$
\aut(\D(\A))=\prod_{C\in\pi_\A}\sym(C)
$$
where $\pi_\A$ is the partition of $\R$ such that $C\in\pi_\A$ iff
$C=\{e\}$ or $C=\{r\in\R^\#:\ d_r=d\}$ with $d\in \D(\A)$.

\thrml{f150205a}
Let $\A=(A,\R)$ be a nontrivial amorphic
$C$-algebra. Then $\aut(\A)=\aut(D(\A))$.
\ethrm
\proof The inclusion $\aut(\A)\le\aut(D(\A))$ is obvious. To prove
the converse inclusion it suffices to verify that
$$
c_{r^g,s^g}^{t^g}=c_{r,s}^t,\quad g\in\sym(D(\A)),\quad
r,s,t\in\R.
$$
However, this follows from (\ref{f070205b}) because $d_{u^g}=d_u$
for all $g\in\aut(D(\A))$.\bull

\sbsn\label{l3}
In this subsection given a multiset $D$ of at least 4 positive real
numbers with $1\in D$, and given an integer $\varepsilon\in\{-1,+1\}$
we construct the class $\A_\varepsilon(D)$ of all isomorphic nontrivial
amorphic $C$-algebras. Let us start with the explicit construction.

\sbst{Algebra $\A_\varepsilon(\R,D)$.} Let $\R$ be a finite set
with distinguished element $e$, and let $d:r\mapsto d_r$ be a bijection
from $\R$ to $D$ such that $d_e=1$. Denote by $A$ the linear space
over $\C$ spanned by $\R$ and by $*$ the identity map of $A$.
Given $r,s,t\in\R$ set $c_{r,s}^t$ to be the real number defined
by condition (C3) if $t=e$, by formula (\ref{f070205a}) with
$n=\sum_{r\in\R}d_r$ if $e\not\in\{r,s,t\}$ and $r\ne s\ne t\ne r$
or $r=s\ne t$ or $r=s=t$, and by equalities (\ref{f050205b})
otherwise. Then the linear space $A$ with the multiplication
defined by (\ref{f090205a}) becomes an (possibly non-associative)
algebra over $\C$ in which the set $\R$ is a linear base.

\lmml{f100205a}
The pair $\A_\varepsilon(\R,D)=(A,\R)$ is a
nontrivial amorphic $C$-algebra.
\elmm
\proof Let us check that
the algebra $A$ is associative. To do this we need to verify that
equalities (\ref{f050205a}) hold for all $r,s,t\in\R$. Clearly, it
is true for $e\in\{r,s,t\}$ or $y=e$ or $r=s=t$. Thus it suffices
to check the cases $r\ne s\ne t\ne r$ and $r=s\ne t$, and in each
of them to verify separately all the possibilities for $y$: $y=r$ or
$y=s$ or $y=t$ or $y\not\in\{r,s,t\}$. All these computations are
straightforward but too long to present here. For example, if
$r\ne s\ne t\ne r$ and $y=s$, equality (\ref{f050205a}) is
equivalent to the following relation:
$$
c_{r,s}^rc_{r,t}^s+c_{r,s}^sc_{s,t}^s+c_{r,s}^tc_{t,t}^s=
c_{r,r}^sc_{s,t}^r+c_{r,s}^sc_{s,t}^s+c_{r,t}^sc_{s,t}^t
$$
which can be verified by the explicit substitution in it the
values of $c_{r,s}^t$ given by (\ref{f070205a}) and
(\ref{f050205b}).

Thus the algebra $A$ is associative. It is easy to see that
conditions (C1), (C2) and (C3) are satisfied. Finally, from the
proof of Theorem~\ref{f060205b} (see the first two lines of
(\ref{f070205b})) it follows that equalities (\ref{f050205c}) hold
and hence condition (C4) is also satisfied. Thus
$\A_\varepsilon(\R,D)$ is a $C$-algebra. The third line of
(\ref{f070205b}) shows that the structure constants of $A$
satisfies (\ref{f050205i}). So the $C$-algebra
$\A_\varepsilon(\R,D)$ is amorphic by Theorem~\ref{f050205g}.\bull \\

From Theorem~\ref{f060205b} it follows that two $C$-algebras
$\A_\varepsilon(\R,D)$ and $\A_{\varepsilon'}(\R',D')$ are
isomorphic iff $\varepsilon=\varepsilon'$ and $D=D'$ (as
multisets). Thus denoting by $\A_\varepsilon(D)$ the class
of all $C$-algebras $\A_\varepsilon(\R,D)$ with fixed $D$ and
$\varepsilon$, we come to the following characterization of
nontrivial amorphic $C$-algebras.

\thrml{f090205c}
The class of all nontrivial amorphic $C$-algebras
equals the union of classes $\A_\varepsilon(\D)$ where $D$ runs
over multisets of at least four positive real numbers with
$1\in D$, and $\varepsilon\in\{-1,+1\}$.\bull
\ethrm

We complete the section by giving an example of amorphic
$C$-algebra generalizing affine $C$-algebras.

\sbst{Generalized affine $C$-algebras.} Let
$\A\in\A_\varepsilon(D)$ be a homogeneous $C$-algebra for some
$\varepsilon$ and $D$. Then the number $d_r$ does not depend on
the choice of the basis element $r$ of $\A$ other than $e$.
Denote it by $d$ and set
$n=n(\A)$. We say that $\A$ is a {\it generalized affine}
$C$-algebra if $d=\sqrt{n}+\varepsilon$. In this case from
(\ref{f070205a}) it follows that
\qtnl{f200205b}
c_{r,s}^t=
\css
1,             &\text{if $|\{r,s,t\}|=3$},\\
\varepsilon+1, &\text{if $|\{r,s,t\}|=2$},\\
\varepsilon(3-\sqrt{n})+1, &\text{if $|\{r,s,t\}|=1$}\\
\ecss
\eqtn
for all basis elements $r,s,t$ of $\A$ other than $e$.
By comparing (\ref{f050205d}) and (\ref{f200205b}) with
$(n,\varepsilon)=(q^2,-1)$ we see that any affine $C$-algebra
corresponding to affine plane of order~$q$, is a generalized
affine $C$-algebra. The converse statement is not true because
there are a lot of integers $q$ for which there is no any
affine plane of order~$q$ \cite[p.63]{DEM}

\section{Homogeneous amorphic $C$-algebras}\label{l4}

\sbsn From the point of view of algebraic combinatorics
$C$-algebras with nonnegative integer structure constants are of
special interest. The following result shows that in studying
amorphic $C$-algebras of such a kind one can restrict himself to
the homogeneous case. Since for a homogeneous $C$-algebra
belonging to a class $\A_\varepsilon(D)$ we have
$D=\{1,d,\ldots,d\}$ where the positive real number $d$ arises
$\nu-1$ times with $\nu$ being the dimension of the algebra, we
write $\A_\varepsilon(\nu,d)$ for the class of all $\A_\varepsilon
(\R,D)$ with $|\R| = \nu$.

\thrml{f150502a}
Let $\A=\A_\varepsilon(\R,D)$ be an amorphic
$C$-algebra for some $\varepsilon$, $\R$ and $D$. Suppose that a
positive real number $d$ has the following property:
\qtnl{f190205a}
\{d_r/d:\ r\in\R^\#\}\subset \N.
\eqtn
Then $\A$
is a fusion of a homogeneous amorphic $C$-algebra belonging to the
class $\A_\varepsilon(\nu,d)$ where $\nu=(n-1)/d+1$ with
$n=n(\A)$.
\ethrm
\proof Set $\A'=\A_\varepsilon(\R',\D')$ where
$\R'$ is a set of cardinality $\nu$ with a distinguished element
$e'$ and $\D'=\{1,d,\ldots,d\}$ ($\nu -1$ times). Then from
Theorem~\ref{f090205c} it follows that $\A'$ is a homogeneous
amorphic $C$-algebra. So $\A'\in\A_\varepsilon(\nu,d)$. This
implies that for any partition $\pi'$ of $\R'$ such that
$\{e'\}\in\pi'$ the $\pi'$-fusion of it is an amorphic
$C$-algebra. Take $\pi'=\{C_r:\ r\in\R\}$ where $C_r$ is a subsedt
of $\R'$ of cardinality  $d_r/d$. Then it is easy to see that
$\A'_{\pi'}\in\A_\varepsilon(\D)$. Taking into account that $\A
\in \A_\varepsilon(D)$, we conclude that $\A$ is isomorphic to
$\A'_{\pi'}$ as desired.\bull

\crllr\label{l5}
Any nontrivial amorphic $C$-algebra with integer
structure constants is a fusion of a generalized affine
$C$-algebra.
\ecrllr
\proof Let $\A$ be a nontrivial amorphic
$C$-algebra. Then by Theorem~\ref{f090205c}
$\A\in\A_\varepsilon(D)$ for some $\varepsilon$ and $D$. Suppose
that all the structure constants of $\A$ are integers. Then from
the first line of (\ref{f070205a}) it follows that the numbers
$d=\sqrt{n}+\varepsilon$ and $d_r/d$ are integers for all basis
elements $r\ne e$ of $\A$. By Theorem~\ref{f150502a} this implies
that $\A$ is a fusion of a homogeneous $C$-algebra
$\A'\in\A_\varepsilon(\nu,d)$ where $\nu=(n-1)/d+1$ with
$n=n(\A)$. Since $n(\A')=n(A)$ it follows
that $\A'$ is a generalized affine $C$-algebra.\bull\\

The hypothesis of Theorem~\ref{f150502a} is satisfied whenever the
structure constants of the $C$-algebra $\A$ are rational numbers.
Indeed, as the number $d$ for which condition (\ref{f190205a})
holds one can choose any rational number $1/N$ where $N$ is a
multiple of the least common multiple of the denominators of the
rational numbers $d_r$. Thus, any amorphic $C$-algebra with
rational structure constants is a fusion of infinitely many
nonisomorphic amorphic homogeneous $C$-algebras.

\sbsn From Theorem~\ref{f150205a} it follows that if $\A$ is a
nontrivial homogeneous amorphic $C$-algebra, then
$\aut(\A)\cong\sym(\R^\#)$. On the other hand, if $\A=(A,\R)$ is
an arbitrary $\C$-algebra and $\pi=\orb(G,\R)$ is a partition of
$\R$ to the set of the orbits of a group $G\le\aut(\A)$, then the
$\pi$-fusion of $\A$ is obviously a $C$-algebra. This implies the
following statement.

\thrml{f150502b} A $C$-algebra $\A=(A,\R)$ is homogeneous and
amorphic iff $\aut(\A)\cong\sym(\R^\#)$.\bull \ethrm

From the main result of the paper \cite{IIM91} it follows that the
sufficiency in Theorem~\ref{f150502b} can be strengthen for the
case when the $C$-algebra $\A$ is the adjacency algebra of an
association scheme. Namely, instead of the assumption that
$\aut(\A)\cong\sym(\R^\#)$ one can suppose that the group
$\aut(\A)$ has an elementary abelian 2-subgroup acting regularly
on the set $\R^\#$. It seems to be true that this result does
valid for arbitrary $C$-algebras (see \cite{B92}). It would be
interesting to describe the properties of the
group $\aut(\A)$ providing that the $C$-algebra $\A$ is
homogeneous and amorphic.

\sbsn We recall that a {\it table} algebra is a $C$-algebra with
nonnegative structure constants (see \cite{AB} and the references
in it). In this subsection we deal with the homogeneous amorphic
table algebras. In contrast to the general case this class of
table algebras is selfdual (concerning the duality theory for
$C$-algebras we refer to \cite{BI}). More exactly, the following
statement holds.

\thrml{f160205a}
Let $\A$ be a homogeneous amorphic table algebra
and let $\widehat\A$ be the dual $C$-algebra of $\A$. Then
$\widehat\A$ is also a homogeneous amorphic table algebra.
Moreover, if $\A\in\A_\varepsilon(\nu,d)$ for some
$\varepsilon,\nu,d$, then $\widehat\A\in\A_\varepsilon(\nu,d)$. In
particular, the table algebras $\A$ and $\widehat\A$ are
isomorphic.
\ethrm
\proof Set $\A=(A,\R)$ and
$\widehat\A=(\widehat A,\widehat R)$. Let $G\le\sym(R)$ be a group
acting regularly on $\R^\#$. Then $G\le\aut(\A)$ by
Theorem~\ref{f150502b}. So by
\cite[Proposition~1.1]{IIM91}\footnote{This result is stated for
adjacency algebras of association schemes, but in fact it is true
for arbitrary $C$-algebra} there exists a group $\widehat
G\le\aut(\widehat\A)$ such that $\widehat G\cong G$ and $\widehat
G$ acts regularly on $\widehat R^\#$. When the group $G$ runs over
all the subgroups of $\sym(R)$ acting regularly on $\R^\#$, the
group $\widehat G$ runs over all the subgroups of $\sym(\widehat
R)$ acting regularly on $\widehat\R^\#$. Thus
$\aut(\widehat\A)\cong\sym(\widehat\R^\#)$ and hence by
Theorem~\ref{f150502b} the $C$-algebra $\widehat\A$ is homogeneous
and amorphic. So by Theorem~\ref{f090205c} we have
\qtnl{f160205d}
\A\in\A_{\varepsilon}(\nu,d),\quad
\widehat\A\in\A_{\widehat\varepsilon}(\nu,\widehat d)
\eqtn
where
$\nu=|\R|=|\widehat\R|$, $d=d_r$ for all $r\in\R^\#$, $\widehat
d=d_{\widehat r}$ for all $\widehat r\in\widehat R^\#$ and
$\varepsilon,\widehat\varepsilon\in\{-1,+1\}$. Moreover, from
\cite[(5.23)]{BI} it follows that
\qtnl{f160205b}
d(\nu-1)+1=\sum_{r\in\R}d_r=\sum_{\widehat
r\in\widehat\R}d_{\widehat r}=\widehat d(\nu-1)+1
\eqtn
from which
we conclude that $d=\widehat d$. Denote by $P$ and $Q$ the
first and the second eigenmatrices of the $C$-algebra $\A$ (the
rows and the columns of the matrix $P$ (resp. $Q$) are indexed by
the elements of $\R$ and $\widehat{\R}$ (resp. $\widehat{\R}$ and
$\R$) respectively). Then from \cite[Theorem~1.2]{IIM91} it
follows that there exists a bijection $\R\to\widehat \R,r\to
\widehat r$ such that $P_{r,\widehat s}=P_{s,\widehat r}$ and
$Q_{\widehat r,s}=Q_{\widehat s,r}$ for all $r,s\in\R$. Set
$n=n(\A)=n(\widehat\A)$. Then by Theorems~5.5 and~5.7 and
equality (5.26) of \cite{BI} we have
$$
\frac{Q_{\widehat s,r}}{d_{\widehat s}}=\frac{\overline
P_{r,\widehat s}}{d_r},\quad
c_{r,r}^r=\frac{d_r^2}{n}\sum_{\widehat s\in\widehat\R^\#}
\frac{1}{d_{\widehat s}^2}Q_{\widehat s,r}Q_{\widehat
s,r}\overline Q_{\widehat s,r}, \quad
c_{\widehat s,\widehat
s}^{\widehat s}=\frac{d_{\widehat s^2}}{n}\sum_{r\in\R^\#}
\frac{1}{d_r^2}P_{r,\widehat s}P_{r,\widehat s}\overline
P_{r,\widehat s}
$$
for all $r,s\in\R^\#$, where the bar means the complex conjugation.
Now the above equalities along with the fact that $d_r=d=\widehat
d=d_{\widehat s}$ for all $r,s\in\R^\#$, imply by condition
(C2) that $c_{\widehat r,\widehat r}^{\widehat r}=c_{r,r}^r$ for
all $r\in\R$. So from (\ref{f070205a}) it follows that
$\varepsilon=\widehat\varepsilon$ and we are done by
(\ref{f160205d}).\bull

\rmrk Theorem~\ref{f160205a} shows that given a homogeneous
amorphic table algebra $\A$, the pair of $C$-algebras $\A$ and
$\widehat\A$ are in positive duality in the sense of \cite{EPV}.
\ermrk

\sbsn In this subsection we find conditions for parameters of a
homogeneous amorphic $C$-algebra
$\A=(A,\R)\in\A_{\varepsilon}(\nu,d)$ that guarantee it to be a
table algebra. To do this we rewrite formulas (\ref{f070205a}) in
terms of $d$, $\varepsilon$ and $n=1+(\nu-1)d$ as follows:
\qtnl{f170205a}
c_{r,s}^t=
\css
\frac{d^2}{(\sqrt{n}+\varepsilon)^2}, &\text{if $|\{r,s,t\}|=3$},\\[4pt]
\phantom{-}\frac{d}{2}\left(1+\frac{2d-n+1}{(\sqrt{n}+\varepsilon)^2}\right),
&\text{if $|\{r,s,t\}|=2$},\\[4pt]
\frac{d-n+1}{2}\left(1+\frac{2d-n+1}{(\sqrt{n}+\varepsilon)^2}\right)
+d-1, &\text{if $|\{r,s,t\}|=1$}\\
\ecss
\eqtn
where $r,s,t\in\R^\#$. Then $\A$ is a table algebra
iff the numbers $c_{r,s}^t$ are nonnegative, which is equivalent
to the two following inequalities:
\qtnl{f170205b}
d\ge\frac{n-1-Q}{2}=-(1+\varepsilon\sqrt{n}),
\eqtn
\qtnl{f170205c}
2d^2-3d(n-1-Q)+(n-1)^2-(n+1)Q\ge 0,
\eqtn
where
$Q=(\sqrt{n}+\varepsilon)^2$. Denote by $\theta$ the
right-hand side of (\ref{f170205b}) and by $d^{\pm}$ the roots of
the quadratic polynomial on $d$ in the left-hand side of
(\ref{f170205c}). Then a straightforward computation shows that
$$
d^{\pm}=
\frac{1+\varepsilon\sqrt{n}}{2}\left(-3\varepsilon\pm\sqrt{9+4\varepsilon\sqrt{n}}\right).
$$
and the system of two inequalities (\ref{f170205b}) and
(\ref{f170205c}) has two families of solutions:
$$
\theta\le d\le d^-,\qquad
d\ge\max\{\theta,d^+\}.
$$
These solutions give a necessary and sufficient condition for the
$C$-algebra $\A$ to be an amorphic homogeneous table algebra.


\begin{thebibliography}{99}

\bibitem{AB}
Z. Arad, Harvey I. Blau, {\em On Table Algebras and Applications
to Finite Group Theory}, J. Algebra, {\bf 138} (1991), 137--185.

\bibitem{AFM}
Z.~Arad, E.~Fisman, M.~Muzychuk, {\em Generalized table algebras},
Israel J. Math., {\bf 114} (1999) 29--60.

\bibitem{B92}
E.~Bannai, {\em On a theorem of Ikuta, Ito and Munemasa}, European
J. Combin., {\bf 13} (1992), 1--3.

\bibitem{BI}
E.~Bannai, T.~Ito, {\em Algebraic combinatorics. I},
Benjamin/Cummings, Menlo Park, CA, 1984.

\bibitem{BMW}
L.~D.~Baumert, W.~H.~Mills, R.~L.~Ward, {\em Uniform cyclotomy},
J. Number Theory, {\bf 14} (1982), 67--82.


\bibitem{DX}
J.~A.~Davis, Q.~Xiang, {\em Amorphic association schemes with
negative latin square type graphs}, arXiv:math/0409401, (2004), 13
pp.

\bibitem{DEM} P. Dembowski, {\em Finite Geometrices}, Springer,
Berlin, {1968}.

\bibitem{DIX}
John D. Dixon, Brian Mortimer {\em Permutation groups,}
Springer-Verlag New York, Inc. (1994).

\bibitem{EPV}
S.~A.~Evdokimov, I.~N.~Ponomarenko, A.~M.~Vershik, {\em
$C$-algebras and algebras in Plancherel duality}, Zapiski
Nauchnykh Seminarov POMI, {\bf 240} (1997), 53--66. English
translation: J. Math. Sci., New York, {\bf 96} (1999), 5,
3478--3485.

\bibitem{GK85}
J.~Y.~Gol'fand, M.~H.~Klin {\em Amorphic cellular rings I}, in:
{\it Investigations in Algebraic Theory of Combinatiorial
Objects}, VNIISI, Moscow, Institute for System Studies, 1985,
32--38( in Russian).

\bibitem{IIM91}
T.~Ikuta, T.~Ito, A.~Munemasa, {\em On pseudo-automorphisms and
fusions of an association scheme}, European J. Combin., {\bf 12}
(1991), 317--325.

\bibitem{I85}
A.~V.~Ivanov, {\em Amorphic cellular rings II}, in: {\it
Investigations in Algebraic Theory of Combinatiorial Objects},
VNIISI, Moscow, Institute for System Studies, 1985, 39--49( in
Russian).

\bibitem{IMY91}
T.~Ito, A.~Munemasa, M.~Yamada, {\em Amorphous association schemes
over the Galois rings of characteristic $4$}, European J. Combin.,
{\bf 12} (1991), 513--526.

\end{thebibliography}
\end{document}